\documentclass{ert-l}

\usepackage[latin1]{inputenc}
\usepackage{amssymb}
\usepackage[all,ps]{xy}
\usepackage{color}

\theoremstyle{plain}
\newtheorem{Def}{Definition}
\newtheorem{theorem}[Def]{Theorem}
\newtheorem{Rem}[Def]{Remark}
\newtheorem{lemma}[Def]{Lemma}
\newtheorem{Prop}[Def]{Proposition}
\newtheorem{Conj}[Def]{Conjecture}
\newtheorem{Not}[Def]{Notation}

\numberwithin{equation}{section}

\newcommand{\rad}{\mathrm{rad}}
\newcommand{\soc}{\mathrm{soc}}

\newcommand{\Hom}{\mathrm{Hom}}
\newcommand{\End}{\mathrm{End}}

\newcommand{\id}{\mathrm{id}}
\newcommand{\im}{\mathrm{im}}
\newcommand{\F}{\mathbb{F}}
\newcommand{\N}{\mathbb{N}}
\newcommand{\Z}{\mathbb{Z}}
\newcommand{\Q}{\mathbb{Q}}

\newcommand{\p}{\mathfrak{p}}
\newcommand{\Oh}{\mathcal{O}}
\newcommand{\A}{\mathcal{A}}
\newcommand{\Ha}{\mathcal{H}}

\newcommand{\idin}{\triangleleft}

\newcommand{\proofbeg}{\noindent{\bf Proof:\ }}
\newcommand{\pend}{\hfill$\blacksquare$}

\begin{document}
\title[Formulas for primitive Idempotents]{Formulas for primitive 
Idempotents in Frobenius Algebras and an Application to Decomposition Maps}
\author{Max Neunhöffer}
\address{School of Mathematics and Statistics, University of St
Andrews, North Haugh, St Andrews, Fife, KY16 9SS, Scotland, United
Kingdom}
\email{neunhoef@mcs.st-and.ac.uk}

\author{Sarah Scherotzke}
\address{Mathematical Institute, 24-29 St Giles' \\ 
Oxford, OX1 3LB, United Kingdom}
\email{scherotz@maths.ox.ac.uk}

\begin{abstract}
In the first part of this paper we present explicit formulas for
primitive idempotents in arbitrary Frobenius algebras using the entries
of representing matrices coming from projective indecomposable modules with
respect to a certain choice of basis. The proofs use a generalisation of
the well known Frobenius-Schur relations for semisimple algebras.

The second part of this paper considers $\Oh$-free $\Oh$-algebras of
finite $\Oh$-rank over a discrete valuation ring $\Oh$ and their
decomposition maps under modular reduction modulo the maximal ideal of $\Oh$,
thereby studying the modular representation theory of such algebras.

Using the formulas from the first part we derive general criteria
for such a decomposition map to be an isomorphism that preserves
the classes of simple modules 
involving explicitly known matrix representations on projective indecomposable
modules.

Finally we show how this approach could eventually be used to attack
a conjecture by Gordon James in the formulation of Meinolf Geck for
Iwahori-Hecke-Algebras, provided the necessary matrix representations
on projective indecomposable modules could be constructed explicitly.
\end{abstract}


\maketitle

\section{Introduction}

Primitive idempotents play a crucial r\^ole in the representation theory
of finite groups and finite-dimensional algebras. In the semisimple case
one has explicit formulas for central primitive idempotents using the
irreducible characters, and for primitive idempotents using the entries
of irreducible matrix representations. The crucial ingredient to prove
these formulas are the well known Frobenius-Schur relations that involve
the matrix coordinate functions of the irreducible representations.

In this paper we generalise this to an arbitrary, finite-dimensional
Frobenius algebra $H$ over a field. To this end, we prove generalisations
of the Frobenius-Schur relations involving the matrix coordinate functions
of representations coming from projective indecomposable modules.
However, we have to choose their basis carefully, namely, the basis must
be adjusted to the socle and the radical of the module, see Section
\ref{averaging} for details.

We then consider algebras over rings and study their decomposition maps
under modular reduction. A general version of Brauer reciprocity shows
that the dual map of such a decomposition map 
can alternatively be defined using idempotents and the corresponding
projective Grothendieck groups. Thus, explicit knowledge about matrix
representations on projective indecomposable modules can
be translated using our formulas for primitive idempotents into knowledge
about decomposition maps.

Finally, we hope that this theory could eventually be applied as part of
a proof of a conjecture by Gordon James in a formulation by Meinolf Geck
about Iwahori-Hecke algebras. 

In Section \ref{averaging} we fix our notation and briefly recall some definitions
and facts. Then we present the well known averaging
operator and some of its properties including the
Theorem by Gasch\"utz and Ikeda for Frobenius algebras.
In the next Section \ref{Rank1} we apply the averaging operator to
linear maps the image of which has dimension $1$ and derive new proofs
in modern language of the generalised Frobenius-Schur relations that
already appear in \cite{Nak} and \cite{Brauer}.

We then proceed to explicitly construct primitive idempotents
using the above relations in Section \ref{primidpot}. An important
feature of these formulas is that we can control the denominators
in the coefficients quite explicitly in terms of the given matrix
representation. This will later be crucial in applying these formulas
to duals of decomposition maps.


Having these preparations at hand, we then set up the concept of
decomposition maps and their duals in Section~\ref{decomp} using the
well-known duality between the Gro\-then\-dieck group $R_0(H)$ of the 
category of finitely generated $H$-modules and the Gro\-then\-dieck
group $K_0(H)$ of the category of finitely generated projective $H$-modules
for a finite-dimensional algebra $H$ over a field. This involves an
extension of the classical
Brauer reciprocity to the non-semisimple case.
In a more general setting, such an extension has
already been discussed in \cite{Centers}. 

The definitions in Section~\ref{appldec} then allow us to use our
formulas for primitive idempotents to derive criteria for a decomposition
map to be trivial, that is, being an isomorphism that preserves the 
classes of simple modules.

Finally, in Section~\ref{james}, we show how this whole theory could 
eventually be applied as a part of a proof for James' conjecture. 
We give another argument that the statement of the conjecture
can only be false for a finite number of cases. However, so far we
can not give an explicit bound. We show in detail what is needed to
apply our results.

\section{Averaging}
\label{averaging}

\medskip
Let $K$ be a field, $H$ a finite-dimensional associative $K$-Algebra and
mod-$H$ the category of finite-dimensional right $H$-modules. We assume
throughout that $K$ is a splitting field for $H$.

For a given $K$-basis $(b_i)_{i=1}^n$ of a right $H$-module $M$
we denote by $(b_i^{\ast})_{i=1}^n$ the basis of the dual space $M^* =
\Hom_K(M,K)$ with $b_i^*(b_j) = \delta_{i,j}$ for $i,j \in \{ 1, \ldots,
n \}$.

\begin{Def}
If $H^*$ contains a linear map $\tau$, such that 
\[ \varphi: H \to H^*, a \mapsto (b \mapsto \tau(ab)) \] 
is a left $H$-module isomorphism, then the
pair $(H,\tau)$ is called a {\bf Frobenius algebra}.
\end{Def}

By \cite[Theorem 61.3]{CR0} the bilinear form $(a|b):=\tau(ab)$ is
non-degenerate and associative if and only if $(H, \tau)$ is a Frobenius
algebra. Therefore we have for every $K$-basis $(C_w)_{w\in W}$ of a
Frobenius algebra $H$ a uniquely determined $K$-basis $(C_y^\vee)_{y \in
W}$ such that $\tau(C_y^{\vee} C_w)=\delta_{y,w}$
for all $y,w\in W$. We call $(C_y^{\vee})_{y\in W}$ the {\bf dual
basis} of $(C_y)_{y\in W}$.

\medskip

In the next lemma we use the following notation. 
Let $\alpha:H \to H$ be an automorphism and $L$
a right $H$-module. Then the vector space $L$ is a right $H$-module with
the following action: $l\ast h:=l\alpha(h)$ for all $l\in L$ and $h\in H$.
We denote the $H$-module $(L,\ast)$ with $L^{\alpha}$.
\begin{lemma}[Frobenius algebras, Nakayama automorphism]
\label{summary}\mbox{}

\begin{enumerate}
\item Let $(H,\tau)$ be a Frobenius algebra, then there exists exactly
one automorphism $\alpha$ of $H$ such that $\tau(ab)=\tau(\alpha(b)a)$ for all
$a,b \in H$. We call $\alpha$ the {\bf Nakayama automorphism}.
\item Let $(H,\tau)$ be a Frobenius algebra, $\alpha$ the Nakayama
automorphism and $P$ a projective indecomposable module of $H$
with socle $S:=\soc(P)$ and head $V:= P/ \rad(P)$. Then $V$ and $S$ are
simple and $V$ is isomorphic to $S^{\alpha^{-1}}$.
\end{enumerate}
\end{lemma}
\proofbeg
\begin{enumerate}
\item As the map $\varphi$ from the first definition is an isomorphism, there is for each
$b \in H$ a unique element $\alpha(b)$ with 
$\tau(ab)=\tau(\alpha(b)a)$ for all $a \in H$. One easily checks that
the map $b \mapsto \alpha(b)$ is linear and we have
$\tau(\alpha(bc)a) = \tau(abc) = \tau(\alpha(c)ab) = \tau(\alpha(b)\alpha(c)a)$
for all $a \in H$. By the non-degeneracy of $\tau$ we conclude
$\alpha(bc) = \alpha(b)\alpha(c)$ for $b,c \in H$ and that $\alpha$
is bijective. Thus $\alpha$ is an automorphism.
\item The modules $V$ and $S$ are both simple by \cite[Prop.~(9.9).(ii)]{CR1}.
Let $e\in H$ be an idempotent with $eH\cong P$. If we identify $P$ and $eH$,
we have $eS=S$. The module $S^{\alpha^{-1}}$ is simple, so
$\Hom_H(eH,S^{\alpha^{-1}})\cong S^{\alpha^{-1}}e$ is
not equal to $\{0\}$ if and only if $S^{\alpha^{-1}}\cong V$. But
$\tau(S^{\alpha^{-1}}e)=\tau(S\alpha^{-1}(e))=\tau(eS)\neq \{0\}$, because
$eS$ is a right ideal and the kernel of $\tau$ does not contain a
non-zero right ideal. It follows that $S^{\alpha^{-1}}e \neq \{0\}$ and
thus $S^{\alpha^{-1}}\cong V$.
\pend
\end{enumerate}

\smallskip

\begin{Not}[Conventions for further reference]
    \label{notation}
In the following we assume that $(H,\tau)$ is a Frobenius algebra,
$\alpha $ the Nakayama automorphism and $P$ a projective indecomposable
module with $\soc(P)=S$ and $P/\rad(P)=V$. Then $P$ is isomorphic to $eH$
for a primitive idempotent $e$. We want to give a formula to compute
such primitive idempotents $e$ in $H$ with $P\cong eH$.

Let $(B_w)_{w\in W}$ be an arbitrary basis of $H$ with dual basis
$(B_w^{\vee})_{w\in W}$. Because of Lemma \ref{summary} part (2)
the dimension $d$ of $V$ is equal to the dimension of $S$. Let $n$ be
the dimension of $P$ and $m := n-d$.
We choose a basis $(b_i)_{i=1}^n$ of $P$ in the following way:

We extend a basis $(b_i)_{i=1}^d$ of $S$ to a basis 
$(b_i)_{i=1}^{m}$ of $\rad(P)$. 
We can then extend the basis $ (b_i)_{i=1}^{m}$
of $\rad(P)$ to a basis $ (b_i)_{i=1}^n$ such that 
$b_i^{\ast}(b_j h)=b_{m+i}^{\ast}(b_{m+j}\alpha(h))$ for all $h \in H$
and $i,j \in \{ 1, \ldots d\}$.
Lemma \ref {summary} part (2) makes sure that there is such a basis.   
\end{Not}

\bigskip

We now introduce the averaging operator, which is similar to the standard
proof of Maschke's Theorem. It will be crucial for our proofs. 

\begin{theorem}[Averaging operator]
\label{averop}
We use the conventions in Notation \ref{notation}.
Let $M,N \in \mathrm{mod}$-$H$ and $f\in \Hom_K (M,N)$. Then the map
$I(f) : M \to N$ with \[I(f)(x):= \sum_{w\in W}f(xB_w)B_w^\vee \mbox{
for all } x\in M\] is a homomorphism of $H$-modules from $M$ to $N$.
Moreover, $I(f)$ does not depend on the choice of basis.

Let $X,Y \in \mathrm{mod}$-$H$ and $\pi \in \Hom_H(X,M)$ and $\psi \in
\Hom_H(N,Y)$ then
\[I(f\circ \pi)=I(f)\circ \pi \mbox{   and   } I(\psi\circ f)=\psi \circ I(f).\] 
\end{theorem}
\proofbeg Straightforward computation. See \cite[Lemma 62.8]{CR0}.
\pend

\begin{theorem}[Gaschütz-Ikeda]
\label{Gachutz-Ikeda}
A right $H$-module $L$ is projective if and only if there is a $\psi\in
\End_K(L)$ with $I(\psi)=\id_L$.
\end{theorem}
\proofbeg \cite[Theorem 62.11]{CR0}
\pend

\begin{lemma}[Averaging homomorphisms between simple modules and PIMs]
\label{1-dim}\mbox{}

We use the conventions in Notation~\ref{notation}.
\begin{enumerate}
\item Let $L$ be a simple right $H$-module which is not isomorphic
to $S$, $f$ a linear form on $L$, $p\in P$ an element of $P$ and
\[\psi:L\to P,\ l\mapsto f(l)p \quad \mbox{ for all } \quad l\in L\] 
 a linear map.  Then $I(\psi)=0$.
 It follows that \[\sum_{w\in W}f(lB_w)p B_w^{\vee}=I(\psi)(l)=0 \quad
  \mbox{ for all }\quad l\in L.\]
\item Let $R$ be a simple right $H$-module which is not isomorphic
to $V$, $f$ a linear form on $P$, $r\in R$ an element of $R$ and
\[\psi:P\to R ,\ p\mapsto f(p)r \quad\mbox{ for all }\quad p\in P\] a
linear map.
 Then $I(\psi)=0$.
 It follows that \[\sum_{w\in W}f(pB_w)r B_w^{\vee}=I(\psi)(p)=0 \quad
 \mbox{ for all } \quad p\in P.\]
\end{enumerate}
\end{lemma}
\proofbeg
\begin{enumerate}
\item The map $I(\psi)$ is a homomorphism of right $H$-modules from $L$
to $P$. As $L$ is simple, $\ker(I(\psi))=L$ or $I(\psi)$ is injective and
thus $I(\psi)(L)\cong L$ . But $S$ is the only simple submodule of $P$
and not isomorphic to $L$. It follows that $I(\psi)=0$.
\item The map $I(\psi)$ is an element of $\Hom_H(P,R)$. As $R$
is a simple module, we have $\im(I(\psi))=0$ or $I(\psi)$ is surjective and 
thus $P/\ker(I(\psi))\cong R$. But $\rad(P)$ is the only maximal submodule of
$P$, so $V$ is the only simple factor module of $P$. Thus $I(\psi)=0$.
\pend
\end{enumerate}
\begin{theorem}[Identity component of endomorphisms]
\label{uniqc}
For every $\psi\in \End_H(P)$ there is a unique constant $c\in K$, such
that 
\[ \im(\psi-c\cdot \id_P)\subseteq \rad(P). \]
For this $c$ we have $(\psi-c\cdot \id_P)(S)=0$.
\end{theorem}
\proofbeg 
As $\psi(\rad(P)) \subset \rad(P)$, the map $\psi$ induces an
endomorphism $\bar \psi:V\to V$. Since $V$ is simple, we have by
Schur's Lemma, that $\bar \psi$ is multiplication by a scalar $c$.
Therefore $\im(\psi-c\cdot \id_P)\subseteq \rad(P)$ holds. Since $P$ is
finite-dimensional and $\psi-c\cdot \id_P$ is not surjective, it follows
that $\psi-c\cdot \id_P$ has a non-trivial kernel. Since $S$ is the only
simple submodule of $P$, it is contained in the kernel of $ \psi-c\cdot
\id_P$.
\pend

\section{Averaging linear maps of rank $1$}
\label{Rank1}
\smallskip

We keep the conventions in Notation~\ref{notation} from 
Section~\ref{averaging}.
By a ``linear map of rank~$1$'' we mean a linear map the image of which
has dimension~$1$. In particular we are interested in linear maps
that expressed as a matrix with respect to the basis $(b_i)_{i=1}^n$
contain only zeros except in one position where they have a one. 
Note that whenever we express our results about endomorphisms
for the convenience of the reader in matrix terms, we
use row convention in the matrices!

\begin{Def}[Constants $c(i,j,s,t)$]
\label{constants}
Let $f_{s,t} \in \End_K(P)$ with $f_{s,t}(p):=b_s^{\ast}(p)b_t$ for 
$s,t \in \{ 1,\ldots,n\}$ and $p\in P$. Note that with respect to the
basis $(b_i)_{i=1}^n$ this linear map corresponds to a matrix that
contains only zeros except for a one in row $s$ and column $t$.
The $f_{s,t}$ form a basis of $\End_K(P)$.
We define 
\[ c(i,j,s,t):= b_j^{\ast}(I(f_{s,t})(b_i))
   = \sum_{w \in W} b_s^{\ast}(b_i B_w) \cdot b_j^{\ast}(b_t B_w^\vee) . \]
Note that if we express $I(f_{s,t})$ in terms of the basis
$(b_i)_{i=1}^n$ as a matrix, then $c(i,j,s,t)$ is the entry in row $i$
and column $j$.
\end{Def}

It is useful to imagine all occurring endomorphisms as matrices, as
always expressed with respect to the basis $(b_i)_{i=1}^n$ using row
convention. Then the elements of $\End_H(P)$ are lower
block-triangular matrices and can be visualised as
follows:

\[ \left[ \begin{array}{c|c|c}
    c \cdot E_d & \rule{2ex}{0mm} 0 \rule{2ex}{0mm} & 
                  0 \rule[-2.5ex]{0mm}{7ex} \\
  \hline
  * & * & 0 \rule[-2.5ex]{0mm}{7ex} \\
  \hline
  * & * & c \cdot E_d  \rule[-2.5ex]{0mm}{7ex}
\end{array} \right] \]
where $c$ is the constant in $K$ from Theorem~\ref{uniqc} 
and $E_d$ stands for a $(d \times d)$-identity matrix.
For the convenience of the reader we briefly indicate in the following 
results about which regions in what matrices we talk.

\begin{lemma}[Frobenius-Schur relations I]
\label{Lem1}
The following equation holds:
\[ c(i,j,s,t) = \left\{ \begin{array}{c@{\quad\mbox{if}\quad}l@{\qquad}l}
  0 & i \le d \mbox{ and } j > d & (1) \\
  0 & i \le m \mbox{ and } j > m & (2) \\
  \delta_{i,j} c(1,1,s,t) & i,j \le d & (3) \\
  \delta_{i,j}c(1,1,s,t) & i,j > m & (4) 
 \end{array} \right.  \]

\end{lemma}
\noindent Note that this lemma is about the parts of the matrix of
$I(f_{s,t})$ that are above the block diagonal and about the upper left
and lower right corner blocks.

\smallskip
\proofbeg 
Let $c_{s,t}$ be the unique constant (see Theorem \ref{uniqc}) such that
\begin{eqnarray*}
 &(\spadesuit)& \quad \im(I(f_{s,t})-c_{s,t}\cdot\id_P)\subseteq \rad(P) 
    \quad \mbox{ and } \\
 &(\clubsuit)& \quad ((I(f_{s,t})-c_{s,t}\cdot\id_P)(S)=0.
\end{eqnarray*}
 It follows immediately from equation $(\clubsuit)$ for $i\le d$ that
\[ c(i,j,s,t)=b_j^{\ast}(I(f_{s,t})(b_i))=b_j^{\ast}(c_{s,t}\cdot
b_i)=c_{s,t}\cdot \delta_{i,j}. \]
 As $c_{s,t}$ only depends on $s$ and $t$, we have
$c_{s,t}=c(1,1,s,t)=c(i,j,s,t)$ for $i=j\le d.$
This gives the third and the first equation.
We know from $(\spadesuit)$ that
\[ c(i,j,s,t)-b_j^{\ast}(c_{s,t}\cdot\id_P(b_i))=
b_j^{\ast}(I(f_{s,t})-c_{s,t}\cdot\id_P)(b_i))=0 \]
if $j>m$.
Thus $c_{s,t}=c(i,j,s,t)=\delta_{i,j}\cdot c_{s,t}$ if $j>m$. As
$c_{s,t}$ only depends on $s$ and $t$ this gives us the last and the
second equation.
\pend

\begin{lemma}[Shifting]
\label{Lem2}
If $s,i\leq d$, then we have $c(i,j,s,t)=c(t,s+m,j,i+m)$.  
\end{lemma}
\noindent Note that this lemma relates certain entries of the matrices
of $I(f_{s,t})$ and $I(f_{j,i+m})$ to each other. Namely, it says that 
for $s \le d$ the entries in the top-most block row of $I(f_{s,t})$
can be determined by looking at the $(t,s+m)$-entry of all the matrices
$I(f_{j,i+m})$ for $1 \le i \le d$ and $1 \le j \le n$.

\smallskip
\proofbeg If $s,i\leq d$, then we have 
\begin{eqnarray*} c(i,j,s,t)&=& 
\sum_{w\in W}b_s^{\ast}(b_iB_w)b_j^{\ast}(b_t B_w^{\vee})
=\sum_{w\in W}b_{s+m}^{\ast}(b_{i+m}\alpha(B_w))b_j^{\ast}(b_t B_w^{\vee})\\
&=&\sum_{w\in W}b_j^{\ast}(b_t B_w^{\vee})b_{s+m}^{\ast}(b_{i+m}\alpha(B_w))
=c(t,s+m,j,i+m),
\end{eqnarray*}
using in the last step the fact that $(\alpha(B_w))_{w\in W}$ is the
dual basis of $(B_w^{\vee})_{w\in W}$ and that the formula for $I$ is
independent of the choice of basis (see Theorem \ref{averop}).
\pend

\begin{lemma}[Frobenius-Schur relations II]
\label{lem3}
We have the following relations:
\[ c(i,j,s,t) = \left\{ \begin{array}{c@{\quad \mbox{if }}l@{\quad}l}
  0 & t\leq m \mbox{ and } ( i,j \le d \mbox{ or } i,j>m ) & (5) \\
  0 & s > d \mbox{ and } ( i,j \le d \mbox{ or } i,j>m )& (6) \\
 \delta_{s+m,t}\delta_{i,j}c & s \leq d \mbox{ and } t>m \mbox{ and }
 i,j\le d& (7) \\
 \end{array} \right.  \]
where we abbreviate $c(1,1,1,1+m)$ to $c$.

\smallskip

Additionally, if $S$ is not isomorphic to $V$ we have:\[ c(i,j,s,t) = \left\{ \begin{array}{c@{\quad\mbox{for}\quad}l@{\qquad}l}
 0 & t \le d & (8)  \\
 0 & s > m &(9) \\
 \end{array} \right.  \]

\end{lemma}
\noindent Note that this lemma is concerned with the upper left and
lower right corner blocks of the matrix of $I(f_{s,t})$ for different
cases of $s$ and $t$. It shows that the result is non-zero only if $s
\le d$ and $t = s+m$, that is, if $f_{s,t}$ has its non-zero entry on 
the diagonal of the upper right corner block.

\smallskip
\proofbeg 
If $t\le m $ and $i,j\le d$ or $i,j>m$, then by Lemma \ref{Lem1}
we have $c(i,j,s,t)=\delta_{i,j}c(1,1,s,t) = \delta_{i,j} c_{s,t}$ where the
constant $c_{s,t}$ is the same as in the proof of Lemma~\ref{Lem1}. Since
$\im(f_{s,t}) \subseteq \rad(P)$ and thus $\im(I(f_{s,t}))\subseteq \rad(P)$, 
the constant $c_{s,t} = c(1,1,s,t)$ is equal to $0$, which proves
$(5)$.

If $s>d$ then $f_{s,t}(S)=0$ and thus $I(f_{s,t})(S)=0$. The map is not
injective and therefore not surjective. So $\im(I(f_{s,t}))\subseteq
\rad(P)$ and analogously $0=c_{s,t}\delta_{i,j}$ for $i,j\le d$ or
$i,j>m$ proving equation $(6)$.

If $s\le d$ and $t>m$ and $i,j\le d$, then we have:
\begin{eqnarray*}
c(i,j,s,t)&=&\delta_{i,j}c(1,1,s,t)=\delta_{i,j}c(t,s+m,1,1+m) \\
 &=&\delta_{i,j}\delta_{s+m,t}c(1,1,1,1+m) 
\end{eqnarray*}
by Lemmas \ref{Lem1} and \ref{Lem2}. 
This proves the statement in equation $(7)$.

\smallskip

If $V$ is not isomorphic to $S$ the set $\{\psi\in
\End_H(P)|\im(\psi)\subseteq S\} = \{0\}$. But
$\im(I(f_{s,t}))\subseteq S$ for $t\le d$ and therefore $I(f_{s,t})=0$.
This gives equation $(8)$.

Furthermore, if $s > m$, then $f_{s,t}$ and thus $I(f_{s,t})$ contain
$\rad(P)$ in their kernel. Therefore the image of $I(f_{s,t})$ is either
isomorphic to the simple module $V$ or is equal to zero. Since the socle
of $P$ is the simple module $S \not\cong V$ it follows that $I(f_{s,t})=0$
and thus $c(i,j,s,t) = 0$ in that case proving equation $(9)$.
\pend
\begin{lemma}
    \label{c}\mbox{}

The constant $c(1,1,1,1+m)$ is not equal to zero.
\end{lemma}
\proofbeg By Theorem \ref{Gachutz-Ikeda}, we know that there is a
linear map $\psi:P\to P$ such that $I(\psi)=\id_P$. It follows that
$b_i^{\ast}(I(\psi)(b_i))=1$ for all $i\in\{1,\ldots,n\}$. The
map $\psi$ can be written as a linear combination of the elements
$(f_{s,t})_{s,t=1}^n$ with coefficients $(d_{s,t})_{s,t=1}^n$. If we
choose $i\le d$ we get $1=\sum_{t=1+m}^n d_{t-m,t}c(1,1,1,1+m)$ using
Lemma \ref{lem3}. This gives $c(1,1,1,1+m)\neq 0$.
\pend

\smallskip
\begin{Rem}\label{noncanonical}
    If we just change the basis $(b_i)_{i=1}^n$ in the
    way that we multiply all basis vectors for $i > m$ with the same 
    non-zero constant then $c(1,1,1,1+m)$ also changes by the same
    constant. At the same time this change does not violate our
    hypotheses in Notation~\ref{notation}. This argument shows that
    even if the constants $c(1,1,1,1+m)$ (for the various projective
    indecomposable modules $P$) look like the Schur elements in the
    semisimple case, they are by no means canonical or invariants of
    the algebra, but depend on the actual choice of the bases in the 
    projective indecomposable modules. For those projective
    indecomposable modules that are simple this constant is the
    Schur element.
\end{Rem}

\section{Primitive Idempotents}
\label{primidpot}

We have now finished all preparations to derive explicit formulas for
primitive idempotents. We continue to use our notation from Sections
\ref{averaging} to \ref{Rank1}. 

\begin{theorem}[Primitive idempotents]
\label{primids}\mbox{}

\begin{enumerate}
\item
Let $c:=c(1,1,1,1+m)$. Then by Lemma \ref{c} we know that $c\neq 0$ and
we can set:
\begin{equation}
\tilde e_i := c^{-1}  \sum_{w \in W} b_i^*( b_{i+m} B_w) 
                                    B_w^{\vee}
\end{equation} for $i\in\{1,\ldots,d\}$. Then there is a polynomial $f
\in \Z[X]$ such that for $e_i := f(\tilde e_i)$ the following is true:
\[
e_i^2 = e_i \quad, \quad  e_i H \cong P \quad \mbox{and} \quad
e_i  e_{i'} \in \rad(H) \]
for $i,i'\in\{1,\ldots,d\}$ and $i \neq i'$.
The polynomial $f$ is given explicitly in the proof.
\item We set \[\tilde E_i := c^{-1}  \sum_{w \in W} b_i^*(
b_{i+m} B_w^{\vee}) B_w \]for $i\in\{1,\ldots,d\}$. Let $\bar
P$ be a projective indecomposable module, the head of which is isomorphic to
$S$.
Then for the $E_i:=f(\tilde E_i)$ the following is true:
\[E_i^2 = E_i \quad, \quad  E_i H \cong \bar P \quad \mbox{and} \quad
E_i  E_{i'} \in \rad(H) \]
for $i,i'\in\{1,\ldots,d\}$ and $i \neq i'$.
\end{enumerate}
\end{theorem}
\noindent Note that the entries of the representing matrices of the
$B_w$ we need are those in the lower left block corner!

\smallskip
\proofbeg \begin{enumerate} \item For $c\neq 0$ see Lemma \ref{c}. In order
to determine the action of $\tilde e_i$ on $V $, we choose $s,j >m$ and
consider:
\begin{eqnarray*} b_s^*\big(b_j  \tilde e_i\big)
   &=& c^{-1}  b_s^*\left( b_j  
     \left( \sum_{w \in W} b_i^*\big( b_{i+m}  B_w \big) B_w^{\vee}
     \right) \right)\\
   &=& c^{-1}  \sum_{w \in W} b_i^*\big(b_{i+m}  B_w\big) 
            b_s^*\big( b_j  B_w^{\vee}\big)\\
&=&c^{-1}c(i+m,s,i,j)=c^{-1}c(j-m,i,s-m,i+m) \\
&=&\delta_{i,j-m}\delta_{s,i+m},
\end{eqnarray*}
using Lemma \ref{Lem2} read backwards and equation $(7)$ of Lemma
\ref{lem3} in the last step. This means that $\tilde e_i$ acts on $V$ as
the projection onto $\left <b_{i+m}+\rad(P) \right >_K$ and thus
$V \tilde e_i$ is one-dimensional. This means that the representing
matrix of $\tilde e_i$ on $V$ with respect to the basis 
\[ (b_{i+m}+\rad(P))_{1 \le i \le d} \] 
is a matrix containing one single
$1$ on the diagonal and apart from that only zeros. So the product
$\tilde e_i  \tilde e_{i'}$ annihilates $V$ for $i\neq i'$.

Let $R$ be a simple module that is not isomorphic to $V$ and $r\in R$
an arbitrary element, then we know by Lemma \ref{1-dim} part (2) with
$f:=b_i^{\ast}$ that
\[ r  \tilde e_i
   = c^{-1}  \sum_{w \in W} b_i^*\big( b_{i+m}  B_w \big)
             r  B_w^{\vee} = 0. \]
Thus $\tilde e_i$ annihilates every simple right module which is not
isomorphic to $V$. Since the elements $\tilde e_i^2 - \tilde e_i$
for $i\in\{1,\ldots,d\}$ and $\tilde e_i  \tilde e_{i'}$ for $
i,i' \in \{1,\ldots ,d\}$ and $i \neq i'$ annihilate every simple right
$H$-module, they are contained in the radical of $H$.

The associative algebra $H$ is  finite-dimensional and thus Artinian. So the radical of $H$ is nilpotent and we can choose an integer $a$ so that 
$(\rad(H))^a=0$. Set 
\[ f := \sum_{s = 0}^a \binom{2a}{s}
   X^{2a-s} (1-X)^s
   \in \Z[X], \]
where $X$ is an indeterminate.
Then with \cite[(6.7)]{CR1} we get the following properties of $f$:
\[ \begin{array}{cl}
\mbox{(i)} & f \in \Z[X] \\
\mbox{(ii)} & f^2 - f \equiv 0 \pmod{X^a(1-X)^a} \\
\mbox{(iii)} & f \equiv X \pmod{X(1-X)}
\end{array} \]
Because of (iii) we know that $e_i := f(\tilde e_i)$ acts on every
simple right $H$-module exactly like $\tilde e_i$. Therefore $e_i \neq 0$ and
$e_i e_{i'}\in \rad(H)$ for $i\neq i'$.
With the choice of $a$ we have 
\[\tilde e_i^a(1-\tilde e_i^a)=(\tilde e_i-\tilde e_i^2)^a=0 \] 
because polynomials in $e_i$ commute. Therefore and because of (ii) the
$e_i$ are idempotents for $i\in\{1,\ldots,d\}$.
This shows that $e_i H$ is a projective $H$-module. Proposition
\cite[(6.6)]{CR1} states that the isomorphism type of a projective
module is determined by the isomorphism type of its head. So $e_i H
\cong P$ if and only if $e_i H /\rad(e_i H) \cong V$.
Let $L$ be an arbitrary simple right $H$-module, then $\Hom_H(e_i H,L)
\cong L e_i$. As $K$ is a splitting field and $e_i H$ is projective, the
$K$-dimension of $\Hom_H(e_i H,L)$ is equal to the multiplicity of $L$
in the head of $e_iH$. We have already shown that $L e_i$ is
one-dimensional for $L\cong V$ and equal to $0$ otherwise. Thus $e_i H$
is indecomposable and isomorphic to $P$.
\item We first determine the action of the $\tilde E_i$ on $S$.
Therefore we choose $j,s\leq d$ and consider:
\begin{eqnarray*}
b_s^*\big(b_j  \tilde E_i\big)
   &=& c^{-1}  b_s^*\left( b_j  
     \left( \sum_{w \in W} b_i^*\big( b_{i+m}  B_w^{\vee}\big) B_w
     \right) \right)\\
   &=& c^{-1}  \sum_{w \in W} b_i^*\big(b_{i+m}  B_w^{\vee}\big)  b_s^*\big( b_j  B_w\big)\\
&=&c^{-1}c(j,i,s,i+m)=\delta_{j,i}\delta_{s,i},
\end{eqnarray*}
using Lemma \ref{lem3} in the last step. So $S\tilde E_i$ is one-dimensional.

Let $L$ be a simple right $H$-module that is not isomorphic to $S$ and
$l\in L$ an arbitrary element. Let $f$ be an arbitrary linear form on $L$,
then
\begin{eqnarray*} f(l  \tilde E_i)
   &=& c^{-1}  \sum_{w \in W} b_i^*\big( b_{i+m}  B_w^{\vee} \big)
             f(l  B_w) \\&=& c^{-1}b_i^{\ast}\left (\sum_{w \in
W}f(l  B_w) b_{i+m}  B_w^{\vee}\right )=0, 
\end{eqnarray*}
using in the last step that $\sum_{w \in W}f(l  B_w) b_{i+m}
 B_w^{\vee}=0$ by Lemma \ref{1-dim} part (1) with $p:=b_{i+m}$.
Since $f$ is an arbitrary linear form we can deduce that $l \tilde
E_i=0$ for all $i\in\{1,\ldots ,d\}$. Thus $\tilde E_i$ annihilates every
simple right module which is not isomorphic to $S$.
The rest of the proof is analogous to the proof in part (1) if we replace
$V$ by $S$.
\pend
\end{enumerate}

\section{Decomposition maps and their duals}
\label{decomp}

This section briefly recalls some notation and definitions needed in
the following sections. All these concepts can be defined in a more
general way, but we do not need them in full generality and thus can
avoid some additional complications, especially due to the fact that we
always assume our base fields to be splitting fields.
We start by introducing some notation for decomposition maps.

For a finite-dimensional algebra $\A$ over a field we denote by
$R_0(\A)$ the Gro\-then\-dieck group of the category of finite-dimensional
right $\A$-modules and by $K_0(\A)$ the Grothendieck group of the
category of finite-dimensional projective right $\A$-modules. For a
right $\A$-module $M$ we denote by $[M]$ its class in $R_0(\A)$ or
$K_0(\A)$, depending on the context. Recall that there is a bilinear
form $\left< - | - \right> : K_0(\A) \times R_0(\A) \to \Z$ given
by $\left< [P] | [V] \right> = \dim \Hom_\A(P,V)$. Then the
set of classes of projective indecomposable $\A$-modules
is the dual basis of the set of classes of simple $\A$-modules
with respect to that form, because for every projective indecomposable
$\A$-module $P$ the head $P/\rad(P)$ is isomorphic to a simple module
and every isomorphism type of simple modules arises in this way
(see \cite[(6.9)]{CR1}).

Let $K$ be a number field, that is, a finite extension of the field
$\Q$ of rational numbers, and $R$ its ring of integers, thus $R$ is a
Dedekind domain and $K$ is the field of fractions
of $R$. For every prime ideal $\p \idin R$ the localisation $R_\p$ of
$R$ at $\p$ is a discrete valuation ring (see \cite[11.2]{Matsu}) and
the residue class field $k_\p := R/\p \cong R_\p/\p_\p$ 
(see \cite[(4.1)]{CR1}) is a finite field of characteristic $\ell$ for
$\ell$ being the rational prime contained in $\p$.

Let $H$ be an associative $R$-free $R$-algebra of finite $R$-rank, and
let $KH := K \otimes_R H$ be its extension of scalars.
Note that we do not assume $KH$ to be semisimple and
in fact the examples occurring in Section~\ref{james} will not be
semisimple. However, we assume that $K$ is a splitting field for $KH$
and that $k_\p$ is a splitting field for the modular reduction
$k_\p H := k_\p \otimes_R H$.

In this situation, the natural map $R_\p \to R_\p/\p_\p \cong k_\p$
gives rise to a ($\Z$-linear) decomposition map $d_\p : R_0(KH) \to
R_0(k_\p H)$ in the following way: A class $[V] \in R_0(KH)$ of a simple
$\A$-module is mapped to $[k_\p \otimes_{R_\p} \tilde V] \in R_0(k_\p
H)$ where $\tilde V$ is an $R_\p H$-module with $K \otimes_{R_\p} \tilde
V \cong V$ as $KH$-modules. Such a module $\tilde V$ exists because $R_\p$ is
a valuation ring and thus every finitely generated torsion-free
$R_\p$-module is free.
See \cite[7.4]{GP} for details on why $d_\p$ is well-defined in this
way. 

\smallskip
We now define a linear map $e_\p : K_0(k_\p H) \to K_0(K H)$,
which will turn out to be closely related to $d_\p$ in the sequel:

\begin{Def}[The dual map of the decomposition map]
\label{map_e}
Let $P = fk_\p H$ be a projective indecomposable $k_\p H$-module
where $f \in k_\p H$ is a primitive idempotent. By \cite[Satz
3.4.1]{Mueller} or \cite[Exercise 6.16]{CR1}
there is an idempotent $e \in R_\p H$ which is mapped to $f$ by the
map $1_{k_\p} \otimes_{R_\p} - : R_\p H \to k_\p H$. We set
$e_\p( [P] ) := [eKH] \in K_0(KH)$.
By standard arguments one shows that this is well-defined
\end{Def}

\smallskip

\noindent The relation between $d_\p$ and $e_\p$ is described by the following
proposition:

\begin{Prop}[Brauer reciprocity]
\label{BrauerRez}
The map \[ e_\p : K_0(k_\p H) \to K_0(KH)\] is the dual map of 
the map \[ d_\p : R_0(KH) \to R_0(k_\p H) \] with respect to the
pairing $\left< - | - \right>$ between $K_0(KH)$ and $R_0(KH)$,
and $K_0(k_\p H)$ and $R_0(k_\p H)$ respectively. More precisely, we have
\begin{equation}
    \label{BrRezEq}
 \left< e_{k_\p}([P]), [V] \right> = \left< [P], d_{k_\p}( [V] ) \right>
\end{equation}
for all $[P] \in K_0(k_\p H)$ and all $[V] \in R_0(KH)$.
\end{Prop}
\proofbeg This is proved in exactly
the same way as \cite[Theorem~18.9]{CR1}. Note that, in this
reference, the algebra $KH$ is globally assumed to be semisimple
(or even a group algebra), but this is completely irrelevant for the
proof of Equation~\ref{BrRezEq}.\pend

\smallskip
Taking the classes of simple modules as bases for
$R_0(KH)$ and $R_0(k_\p H)$ we can express the decomposition map
$d_\p$ as a matrix, the so called ``decomposition matrix of $d_\p$'',
the rows of which are indexed by the basis of $R_0(KH)$ and the columns
of which are indexed by the basis of $R_0(k_\p H)$. A row of the
decomposition matrix thus contains the multiplicities of the simple
$k_\p$-modules in a modular reduction of the corresponding simple
$KH$-module.

Analogously, if we take the classes of projective indecomposable modules
as bases for $K_0(k_\p H)$ and $K_0(KH)$, we can express the map $e_\p$
as a matrix, the rows of which are indexed by the basis of $K_0(k_\p
H)$ and the columns of which are indexed by the basis of $K_0(KH)$. A
row of this matrix thus contains the multiplicities of the projective
indecomposable $KH$-modules in a direct sum decomposition of a lift of a
projective indecomposable $k_\p H$-module.

Since our chosen bases of $K_0(KH)$ and $K_0(k_\p H)$ are just the
dual bases of our chosen bases of $R_0(KH)$ and $R_0(k_\p H)$ with
respect to the pairing $\left< - | - \right>$, Proposition~\ref{BrauerRez}
states that the matrix of $e_\p$ is just the transposed matrix of
the decomposition matrix. Thus a column of the decomposition matrix
contains the multiplicities of the projective indecomposable $KH$-modules
in a direct sum decomposition of a lift of a projective indecomposable
$k_\p H$-module, which is nothing but the classical Brauer reciprocity.
See \cite[7.5.2]{GP} for the corresponding result if
$KH$ is semisimple, and \cite[Section 2]{Centers} for a more general
and more complicated case not assuming splitting fields. Our exposition
in the present paper is similar to \cite[Kapitel V]{PhdMax}.

\section{Application to decomposition maps}
\label{appldec}

This section is motivated by James' conjecture on Iwahori-Hecke algebras.
The setup presented here is a generalisation of the situation in 
Section~\ref{james}.

Let $R$ be a Dedekind domain, $K$ its field of fractions, $H$ an
associative $R$-free $R$-algebra with finite $R$-rank such that its
extension of scalars $KH := K \otimes_R H$ is a Frobenius algebra with
$K$-linear map $\tau : KH \to K$. Note that we do not assume $KH$ to be
semisimple and in fact the
examples occurring in Section~\ref{james} will not be semisimple.

We first formulate a criterion for a column in the decomposition matrix
to be trivial:

\begin{lemma}[Trivial columns in the decomposition matrix]
\label{crit}
Let $\p$ be a prime ideal of $R$ and $e \in KH$ a primitive idempotent
that lies in $R_\p H := R_\p \otimes_R H$. Then the idempotent
$f := 1_{k_\p} \otimes_{R_\p} e \in k_\p H$ is primitive and
the column of the decomposition matrix of $d_\p$ 
(see the end of Section~\ref{decomp})
corresponding to the simple $k_\p H$-module $f k_\p H / \rad( f k_\p H)$
contains exactly one $1$ in the row corresponding to the simple $KH$-module
$e KH / \rad( e KH )$ and apart from that only zeros.
\end{lemma}
\proofbeg
We use Brauer reciprocity as described in Section~\ref{decomp}: Since
$R_\p H$ is semiperfect by 
\cite[Satz~4.3.1]{Mueller} or \cite[Exercise~6.16]{CR1}
the idempotent
$f := 1_{k_\p} \otimes_{R_\p} e \in k_\p H$ is primitive
and the class $[fk_\p H] \in K_0(k_\p H)$ is mapped to 
the class $[e KH] \in K_0(KH)$ of the projective indecomposable
module $eKH$ by the map $e_\p$. The interpretation of $e_\p$ as the
dual map of the decomposition map $d_\p$ and the definition of the
decomposition matrix of $d_\p$ gives the statement in the lemma.
\pend

\medskip
Now we use a concrete representation of $KH$ on a projective indecomposable
module $P$ to find an infinite set of prime ideals $\p \idin R$ for which
Lemma~\ref{crit} can be applied.

\smallskip
Let $(B_w)_{w \in W}$ be an $R$-basis of $H$. Then it is also a $K$-basis
of $KH$ and we denote its dual basis with respect to $\tau$ by
$(B_w^\vee)_{w \in W}$.

Choose a $K$-basis $(b_1,\ldots,b_d,\ldots,b_m, \ldots, b_n)$ of $P$ as in 
Section~\ref{averaging}. Scale the basis vectors $(b_{m+1}, \ldots, b_n)$ 
corresponding to the head of $P$ by a common scalar in $K$, such that the 
constant $c(1,1,1,1+m)$ (see Definition~\ref{constants} and
Remark~\ref{noncanonical}) is equal to
$1$. For a fixed $1 \le i \le d$, write every number 
$a_w := b_i^*(b_{i+m} \cdot B_w)$ for $w \in
W$ as a quotient $a_w =: s_w / t_w$ with $s_w \in R$ and $t_w \in R
\setminus\{0\}$ and let $I_i \idin R$ be the ideal generated by the
product $\prod_{w \in W} t_w$ of all such denominators. Let $I$ be the
ideal generated by all the $I_i$ for $1 \le i \le d$.

Note that the ideal $I$ depends on the choice of the numerators and
denominators $s_w$ and $t_w$ and of course on the choice of basis
$(b_j)_{1 \le j \le n}$. However, for all such choices, we have:

\begin{theorem}[Criterion for trivial column I]
\label{coltriv1}
Let the ideal $I$ be defined as above
and $\p \idin R$ be a prime ideal such that $I$ is not contained in $\p$.
Then there is a primitive idempotent $e \in R_\p H$ satisfying $eKH \cong P$
as $KH$-modules,
such that the idempotent $f := 1_{k_\p} \otimes_{R_\p} e \in k_\p H$
is also primitive. Furthermore, the column of the decomposition matrix 
of $d_\p$ corresponding to the simple $k_\p H$-module $fk_\p H/\rad(fk_\p H)$
contains only zeros, except in the row corresponding to the simple
$KH$-module $eKH/\rad(eKH)$, where it contains a 1.
\end{theorem}
\proofbeg At least one of the ideals $I_i$
is not contained in $\p$ and thus all the denominators $t_w$ the product
of which generates $I_i$ are not
in $\p$, because $\p$ is a prime ideal. Therefore Theorem~\ref{primids} gives
an idempotent $e_i \in R_\p H$ with $e_i KH \cong P$ as $KH$-modules and 
Lemma \ref{crit} shows the last statement in the theorem.
\pend

\medskip
\textbf{Remark:} As every ideal $I \idin R$ is divided by only finitely
many prime ideals $\p \idin R$, there is only a finite 
number of cases, in which the hypotheses of Theorem~\ref{coltriv1} is
not fulfilled. Repeating this argument for every isomorphism type
of projective indecomposable $KH$-modules shows that there are only
finitely many prime ideals $\p \idin R$ for which the decomposition
matrix is not equal to an identity matrix.

\medskip
We now follow a different approach. Whereas for Theorem~\ref{coltriv1}
we used only some entries of all representing matrices for a basis
$(B_w)_{w \in W}$, we now consider all entries of representing matrices
for a system of generators of $H$.

To this end, let $X \subseteq W$ be a set, such that $\{ B_x \mid x \in
W\}$ generates $H$ as an algebra. Write every entry $b_i^*(b_j \cdot B_x)$
of the representing matrices for all $1 \le i,j \le n$ and all $x \in X$
as a quotient $s_x/t_x$ and let $N$ be the product of all the denominators
$t_x$. Let $J := NR \idin R$ be the ideal of $R$ generated by $N$.

Note that the ideal $J$ depends on the choice of the numerators and
denominators, of course on the choice of the basis $(b_j)_{1 \le j \le
n}$, and on the choice of the generating system $\{B_x \mid x \in X\}$. 
However, for all such choices we have:

\begin{theorem}[Criterion for trivial column II]
\label{coltriv2}
Let the ideal $J$ be defined as above
and $\p \idin R$ be a prime ideal such that $J$ is not contained in $\p$.
Then there is a primitive idempotent $e \in R_\p H$ satisfying $eKH \cong P$
as $KH$-modules,
such that the idempotent $f := 1_{k_\p} \otimes_{R_\p} e \in k_\p H$
is also primitive. Furthermore, the column of the decomposition matrix 
of $d_\p$ corresponding to the simple $k_\p H$-module $fk_\p H/\rad(fk_\p H)$
contains only zeros, except in the row corresponding to the simple
$KH$-module $eKH/\rad(eKH)$, where it contains a 1.
\end{theorem}
\proofbeg All denominators of all entries
of the representing matrices of all generators $B_x$ for $x \in X$
are not in $\p$ and so these matrix entries lie in $R_\p$. Since
the $B_x$ generate $H$ as an algebra, the same holds for all
representing matrices of all elements $B_w$ for $w \in W$.
Thus Theorem~\ref{primids} gives
an idempotent $e_1 \in R_\p H$ with $e_1 KH \cong P$ as $KH$-modules and 
Lemma \ref{crit} shows the last statement in the theorem.
\pend

\medskip
\textbf{Remark:} The same argument as the one after Theorem~\ref{coltriv1}
applies here, showing that the hypotheses of Theorem~\ref{coltriv2}
are fulfilled for all but a finite number of prime ideals $\p$.

%
%
%
%
%
%

\section{A possible application to James' conjecture}
\label{james}

In this section we present the situation in which the above results
could be applied, provided one could construct certain
representations, such that one could control the denominators
of the entries in the representing matrices.

Let $(W,S)$ be a finite Coxeter system, that is, $W$ is a finite group
with a subset $S$ such that we have a presentation of the form
\[ W = \left< s \in S \mid s^2 = 1 \mbox{ and }
   (st)^{m_{s,t}} = 1 \mbox{ for } s,t \in S \right>, \]
where $m_{s,t}$ is the order of $st$.
Let $L : W \to \N \cup \{0\}$ be the length function on $W$, that is,
$L(w)$ is the number of factors in the shortest expression of $w$
as product of generators in $S$. An expression 
$w = s_1 \cdot \cdots \cdot s_{L(w)}$ with $s_i \in S$ is called
\emph{reduced}.

Let $A$ be any commutative ring with $1$ and $v \in A$ invertible.
We can now define the one-parameter 
Iwahori-Hecke algebra $\Ha_A(W,S,v)$
over $A$ to be the associative
$A$-algebra with generators $\{T_w \mid w \in W\}$ subject to the relations
\[ \begin{array}{rcl@{\hspace*{5mm}}l}
  T_s^2 &=& T_\id + (v-v^{-1}) T_s & \mbox{for all } s \in S \\
   T_w &=& T_{s_1} \cdot \, \cdots \, \cdot T_{s_k}
       & \mbox{for every reduced expression } \\
         &&&w = s_1 \cdot\,\cdots\,\cdot s_k \mbox{ in } W 
        \mbox{ with } s_i \in S 
   \end{array} \]
where $\id \in W$ denotes the identity element.

The algebra $\Ha$ is free as an $A$-module with basis $(T_w)_{w \in W}$ 
(see \cite[3.3]{Uneq}) and has a symmetrising trace map $\tau: \Ha \to A,
T_\id \mapsto 1, T_w \mapsto 0$ for $\id \neq w \in W$, 
which makes $\Ha$ into a symmetric algebra in the sense of
\cite[7.1.1]{GP}.
The dual basis of $(T_w)_{w \in W}$ with respect to $\tau$ is
$(T^\vee_w)_{w \in W}$ with $T^\vee_w = T_{w^{-1}}$ (for all of
this, see \cite[10.3,10.4]{Uneq}).
Note that in Section~\ref{averaging} we chose to define Frobenius algebras
only over fields. If $A$ is a field, $(\Ha,\tau)$ is just a Frobenius
algebra with the additional property that $\tau(xy) = \tau(yx)$ for all
$x,y \in \Ha$.

This construction is functorial in the sense, that if $f : A \to B$ is
a homomorphism into a commutative ring $B$, we can regard $B$ as an
$A$-module via $f$ and then have a canonical isomorphism
$B \otimes_A \Ha_A(W,S,v) \cong \Ha_B(W,S,f(v))$, where the
latter is defined in exactly the same way as above as a finitely
presented algebra, only over the ring $B$ with parameter $f(v)$ instead
over $A$ with parameter $v$ (see \cite[8.1.2]{GP}).

We now consider three different base rings:
Firstly let $\tilde K \supseteq \Q$ be a finite extension such that
$\tilde K$ is a splitting field for the group algebra $\tilde K W$. 
Let $\tilde R$ be the ring of integers of $\tilde K$, it is a Dedekind
domain and in particular integrally closed. Set $\hat R := \tilde R[v,v^{-1}]$,
the ring of Laurent polynomials over $\tilde R$. 
Then $\hat K := \tilde K(v)$ is the field of fractions of $\hat R$.
Note that $\hat K$ is a splitting
field for the extension of scalars $\hat K\Ha := \Ha_{\hat K}(W,S,v)
= \hat K \otimes_{\hat R} \Ha_{\hat R}(W,S,v)$.
This follows from \cite[9.3.5]{GP} and the fact that our parameter
$v$ is the square root of the parameter there. In addition, all
irreducible characters of $\hat K\Ha$ can be realized over $\hat K$.

Secondly, we consider a finite field $k$ of characteristic $\ell > 0$ and
a homomorphism $\theta_\ell : \hat R \to k$ such that $k$ is the field of
fractions of the image $\theta_\ell(\hat R)$. Since $v$ is invertible in
$\hat R$, this also holds for the image $q := \theta_\ell(v)$, which thus has
finite order. Let $e$ be the order of $q$ for $q \neq 1$, and else set
$e := \ell$.

Thirdly, let $K$ be $\tilde K(\zeta)$ where $\zeta$ is a primitive
$e$-th root of unity. Let $R$ be the ring of integers of $K$. Then $R$
is a Dedekind domain, which contains $\tilde R$ as a
(possibly equal) subring. Note that the choice of $\zeta$ and thus of
$R$ and $K$ is determined entirely by $k$ and $q$.

In this situation, we have a ring homomorphism 
$\theta_\zeta : \hat R \to R$ mapping $v$ to $\zeta$. Since both $R$ and
$k$ are integral domains, the kernels of $\theta_\ell$ and $\theta_\zeta$
are both prime ideals. Furthermore, the kernel of $\theta_\zeta$ is generated
by the minimal polynomial of $\zeta$ over $\tilde R$, which is an
irreducible factor of the cyclotomic polynomial $\phi_e(v) \in R[v]$
having $\zeta$ as a root. Since $q = \theta_\ell(v)$ has multiplicative
order $e$ (for $q \neq 1$), the field $k$ contains
primitive $e$-th roots of unity and thus
the $\ell$-modular reduction of the cyclotomic polynomial 
$\phi_e \in \Z[X]$ is equal to a product of linear factors over $k$.
Therefore, the kernel of $\theta_\zeta$ is contained in the kernel of
$\theta_\ell$ and there is a ring homomorphism 
$\theta_\ell^\zeta : R \to k$ with 
$\theta_\ell = \theta_\ell^\zeta \circ \theta_\zeta$ thus mapping $\zeta$
to $q = \theta_\ell(v) \in k$. The same holds for the case $e = \ell$
and $q = 1$.

These ring homomorphisms together with the functoriality above gives us
three Iwahori-Hecke algebras $\Ha_{\hat R}(W,S,v)$, $\Ha_{R}(W,S,\zeta)$
and $k \Ha := \Ha_k(W,S,q)$ together with canonical maps between them.
The first two are associative algebras over the rings $\hat R$ and
$R$ respectively such that we can also consider the corresponding
extensions of scalars $\hat K \Ha := \Ha_{\hat K}(W,S,v)$ and $K \Ha :=
\Ha_K(W,S,\zeta)$ to the respective fields of fractions $\hat K$ and
$K$.

For this situation, Meinolf Geck and Rapha\"el Rouquier in \cite[2.5]{Centers}
have defined a commutative diagram of decomposition maps
\[ \xymatrix{ R_0(\hat K \Ha) \ar[rr]^{d_{\theta_\ell}} 
              \ar[dr]_{d_{\theta_\zeta}} && R_0(k\Ha) \\
              & R_0(K\Ha) \ar[ur]_{d_{\theta_\ell^\zeta}} & } \]
We do not want to go into the details of the definition of $d_{\theta_\ell}$
and $d_{\theta_\zeta}$ here, because we do not need those two maps in the
sequel and our definition of decomposition maps in Section~\ref{appldec}
would have to be generalised to do this.

However, we can now formulate Geck's version of James' conjecture
(see \cite[(3.4)]{Bourb}).

\begin{Conj}[{see \cite[Section 4]{James}, \cite[(3.4)]{Bourb}}]
\label{jamesconj}
If in the situation above, $k = \F_\ell$ is the field of $\ell$ elements,
and $\ell$ is coprime to the order of $W$, then 
\smash{$d_{\theta_\ell^\zeta}$} is
an isomorphism that preserves the classes of simple modules.
\end{Conj}

\textbf{Remark 1:} The authors do not see a reason for the restriction
$k = \F_\ell$ in this conjecture. However, there does not seem to be
much computational evidence available showing that the conjecture
holds in cases $k \neq \F_\ell$. Thus we stick to this restriction.

\smallskip
\textbf{Remark 2:} The statement simply means that the
decomposition matrix corresponding to \smash{$d_{\theta_\ell^\zeta}$} is an 
identity matrix. This in turn implies that the modular decomposition map
\smash{$d_{\theta_\ell}$} is completely determined by the
decomposition map \smash{$d_{\theta_\zeta}$} which is defined
only involving rings of characteristic $0$.

\smallskip
\textbf{Remark 3:} The definition of $d_{\theta_\ell^\zeta}$ in 
\cite{Bourb} coincides with our definition of a decomposition map
$d_\p : R_0(K\Ha) \to R_0(k\Ha)$ using the kernel of the surjective
ring homomorphism $\theta_\ell^\zeta$ as the prime ideal $\p \idin R$.

\smallskip
To the best knowledge of the authors, this conjecture is still open except
in the defect $1$ case and a few other special cases (see for example 
\cite[6.4]{Chuang}, \cite{BrauerTrees} and \cite{Fayers}).

Using any set $\{P_1, P_2, \ldots, P_t\}$ of representatives of the
isomorphism classes of projective indecomposable $K\Ha$-modules 
together with the argument after the proof of Theorem \ref{coltriv1} we
can show that
there are only finitely many prime ideals $\p \idin R$ and thus only
finitely many pairs $(k,q)$ such that the decomposition map 
$d_\ell^\zeta$ does not fulfil the statement of Conjecture~\ref{jamesconj}.
This reasoning provides an alternative to the corresponding proof in 
\cite[5.5]{BrauerTrees}. 

However, there is no explicit lower bound $B$ for $\ell$ known such that
for all $\ell > B$ and all $(k,q)$ with $\mathrm{char}\ k = \ell$ the
decomposition map $d_\ell^\zeta$ is an isomorphism that preserves
the classes of simple modules.

Our Theorems \ref{coltriv1} and \ref{coltriv2} could eventually be used
to achieve such a bound in the following way: If one could explicitly
construct realisations of projective indecomposable $K\Ha$-modules
together with bases adapted to socle and radical, \textbf{and} one
could control the denominators of the resulting representing matrix
entries (either some specific entries for all representing matrices
of a basis of $\Ha_R(W,S,\zeta)$ (as for Theorem \ref{coltriv1}) or 
for all entries 
for representing matrices of a generating system of $\Ha_R(W,S,\zeta)$ 
(as for Theorem \ref{coltriv2}), then one
would get a criterion for which prime ideals $\p \idin R$ \textbf{all}
columns of the decomposition matrix are trivial and thus the decomposition
map $d_\ell^\zeta$ is an isomorphism that preserves the classes of
simple modules.

It is our hope that such projective modules can be constructed eventually
thus yielding at least a part of a proof of James' conjecture. A
slight hint in this direction is provided by some observations the 
first author has made in his PhD-thesis (see \cite[Section VI.7]{PhdMax}).
There it is reported, that explicit matrix representations coming from
projective modules for the Iwahori-Hecke algebra $\Ha_R(W,S,\zeta)$ over
the ring of integers $R$ of the number field $\Q(\zeta)$ can be obtained
using the Kazhdan-Lusztig basis and intervals in the poset of left
cells and their corresponding non-simple cell modules (see 
\cite[Chapter 8]{Uneq} or \cite{KLWedder} for definitions of these
concepts).

Although these observations are still only the result of
a few computer calculations it seems not totally impossible that the
Kazhdan-Lusztig theory and in particular the methods involving cells
could eventually lead to the explicit construction of projective
modules, in particular since the Kazhdan-Lusztig basis has already been
used extensively to study the representation theory of these algebras.

\bigskip
We kindly thank the referee for his valuable comments and
suggestions which were very useful to streamline the exposition.

\bibliographystyle{alpha}

\end{document}